\documentclass[11pt]{amsart}

\usepackage{amsmath}
\usepackage{amssymb}
\usepackage{color}
\usepackage[colorlinks,citecolor=blue,urlcolor=black,linkcolor=black]{hyperref}

\newtheorem{theorem}{Theorem}[section]
\newtheorem{claim}[theorem]{Claim}
\newtheorem{fact}[theorem]{Fact}

\newtheorem{lemma}[theorem]{Lemma}

\newtheorem{proposition}[theorem]{Proposition}
\newtheorem{corollary}[theorem]{Corollary}

\theoremstyle{definition}
\newtheorem{definition}[theorem]{Definition}

\newtheorem{question}[theorem]{Question}

\theoremstyle{remark}

\newcount\skewfactor
\def\mathunderaccent#1#2 {\let\theaccent#1\skewfactor#2
\mathpalette\putaccentunder}
\def\putaccentunder#1#2{\oalign{$#1#2$\crcr\hidewidth
\vbox to.2ex{\hbox{$#1\skew\skewfactor\theaccent{}$}\vss}\hidewidth}}
\def\name{\mathunderaccent\tilde-3 }

\def\smallbox#1{\leavevmode\thinspace\hbox{\vrule\vtop{\vbox
   {\hrule\kern1pt\hbox{\vphantom{\tt/}\thinspace{\tt#1}\thinspace}}
   \kern1pt\hrule}\vrule}\thinspace}

\newcommand{\cf}{{\rm cf}}

\newcommand{\stick}{\ensuremath\mspace{2mu}\mid \mspace{-12mu}{\raise 0.4em \hbox{$\bullet$}}}

\def\qedref#1{$\qed_{\reforiginal{#1}}$}

\title{Tiltan}
\author{Shimon Garti}
\address{Institute of Mathematics,
 The Hebrew University of Jerusalem,
 Jerusalem 91904, Israel}
\email{shimon.garty@mail.huji.ac.il}
\thanks{}


\begin{document}
\let\labeloriginal\label
\let\reforiginal\ref

\begin{abstract}
We prove the consistency of $\clubsuit$ with the negation of Galvin's property.
On the other hand, we show that superclub implies Galvin's property.
We also prove the consistency of $\clubsuit_{\kappa^+}$ with $\mathfrak{s}_\kappa>\kappa^+$ for a supercompact cardinal $\kappa$.
\end{abstract}

\maketitle

\newpage

\section{Introduction}

The diamond principle of Jensen, \cite{MR0309729}, is a prediction principle. It says that there exists a sequence of sets $(A_\alpha:\alpha\in\omega_1)$, such that each $A_\alpha$ is a subset of $\alpha$, and such that for every $A\subseteq\omega_1$ the set $\{\alpha\in\omega_1:A\cap\alpha=A_\alpha\}$ is a stationary subset of $\omega_1$.

A weaker prediction principle, denoted by $\clubsuit$, was introduced by Ostaszewski in \cite{MR0438292}. Usually it is called the club principle, but we shall employ the name \emph{tiltan} to refer to $\clubsuit$. The reason is that we are going to deal extensively with closed unbounded sets using the acronym \emph{club}, and anticipating a natural confusion we prefer a linguistic distinction.
The name tiltan means clover in Mishnaic Hebrew (but in good old manuscripts it is pronounced taltan, see \cite{Yeivin}, p. 984).
The local version of tiltan at $\aleph_1$ reads as follows:

\begin{definition}
\label{deftiltan} Tiltan. \newline
There exists a sequence of sets $\langle T_\alpha: \alpha\in\lim(\omega_1)\rangle$, such that each $T_\alpha$ is a cofinal subset of $\alpha$, and such that for every unbounded set $A\subseteq\omega_1$ there are stationarily many ordinals $\alpha$ for which $T_\alpha\subseteq A$.
\end{definition}

The definition generalizes easily to any stationary set $S$ of any regular uncountable cardinal $\kappa$ whose elements are limit ordinals. The tiltan sequence will be $\langle T_\alpha: \alpha\in S\rangle$, and the assertion will be denoted by $\clubsuit_S$. Clearly, if $\clubsuit_{S_0}$ holds and $S_1\supseteq S_0$ then $\clubsuit_{S_1}$ holds as well.

It is clear from the definition that $\Diamond \Rightarrow \clubsuit$. The difference is two-fold. Firstly, the $\Diamond$-prediction is accurate and based on the equality relation, namely $A\cap\alpha = A_\alpha$ while the $\clubsuit$-prediction promises only inclusion, i.e. $T_\alpha\subseteq A$. Secondly, the diamond predicts \emph{all} the subsets of $\omega_1$ (or larger cardinals) including the countable subsets, while the tiltan predicts only \emph{unbounded} subsets of $\omega_1$. In both points, $\clubsuit$ is weaker than $\Diamond$.

One may wonder if the tiltan is strictly weaker than the diamond. It is easy to show that $\clubsuit + 2^\omega=\omega_1$ is equivalent to $\Diamond$. The question reduces, therefore, to the possible consistency of tiltan with $2^\omega>\omega_1$. The answer is yes, as proved by Shelah in \cite{MR1623206}, I, \S 7. The proof shows, in particular, the consistency of $\clubsuit + \neg\Diamond$. This result opens a window to a variety of consistency results of this form.

Suppose that $\varphi$ is a mathematical statement which follows from the diamond. One may ask whether the tiltan is consistent with $\neg\varphi$. We focus, in this paper, on a statement which we call \emph{Galvin's property}. It is based on a theorem of Galvin which appears in \cite{MR0369081}. We quote the version of $\aleph_1$ and club sets, but the theorem generalizes to every normal filter over any regular uncountable cardinal.

\begin{theorem}
\label{thmgalvin} CH and Galvin's property. \newline
Assume that $2^\omega = \omega_1$. \newline
Then any collection $\{C_\alpha:\alpha\in\omega_2\}$ of club subsets of $\aleph_1$ admits a sub-collection $\{C_{\alpha_\beta}:\beta\in \omega_1\}$ whose intersection is a club subset of $\aleph_1$.
\end{theorem}

\hfill \qedref{thmgalvin}

Galvin's property follows from CH, and a fortiori from the diamond. Question 2.4 from \cite{MR3604115} is whether Galvin's property follows from the tiltan. The original proof of Galvin gives the impression that the answer should be positive. Surprisingly, we shall prove the opposite by showing the consistency of tiltan with the failure of Galvin's property.

Nevertheless, something from the natural impression still remains and can be proved.
Tiltan is consistent only with a weak negation of Galvin's property. The strong negation of it cannot be true under the tiltan assumption.
Let us try to clarify this point.

Galvin's property deals with a sub-collection whose intersection is a club, but the real point is only unboundedness. If $C = \bigcap\{C_{\alpha_\beta}: \beta\in\omega_1\}$ and $a\subseteq C$ is unbounded, then $c\ell(a)\subseteq C$ as well. Consequently, if one wishes to force the negation of Galvin's property then a bounded intersection must be forced. This is done, twice, in a work of Abraham and Shelah \cite{MR830084}. Our purpose is to combine the forcing of \cite{MR830084} with the classical way to force $\clubsuit + \neg {\rm CH}$, thus obtaining the main result of the next section:

\begin{theorem}
\label{thmmt} It is consistent that $\clubsuit$ holds, $2^\omega=\lambda, \lambda$ is arbitrarily large and there exists a collection $\{C_\alpha: \alpha<\lambda\}$ of club subsets of $\aleph_1$ such that any $\aleph_1$-sub-collection of it has bounded intersection.
\end{theorem}

\hfill \qedref{thmmt}

The negation of Galvin's property, reflected in the above theorem, is different from the situation in \cite{MR830084} notwithstanding. In the constructions of \cite{MR830084} not only any sub-collection of size $\aleph_1$ has bounded intersection (in $\aleph_1$) but it has \emph{finite} intersection. Let us call this property \emph{a strong negation of Galvin's property}. We shall see that tiltan is incompatible with such a strong negation. Namely, under $\clubsuit$ any collection of the form $\{C_\alpha: \alpha<\omega_2\}$ contains even a sub-collection of $\aleph_2$-many sets with infinite intersection. Actually, an intersection of order type $\geq\tau$ for every ordinal $\tau\in\omega\cdot\omega$ can be shown to exist.
This means that the main theorem is optimal in some sense. Moreover, it gives some information about possible ways to force tiltan and their limitations. One way to demonstrate this observation is to strengthen tiltan, as done in the second section. We shall work with the prediction principle \emph{superclub} from \cite{primavesi} and show that it implies Galvin's property. In the last section we deal with the splitting number $\mathfrak{s}_\kappa$ and the possibility that $\clubsuit_{\kappa^+}$ be consistent with $\mathfrak{s}_\kappa>\kappa^+$.

Our notation is mostly standard. If $\kappa=\cf(\kappa)<\lambda$ then $S^\lambda_\kappa = \{\delta\in\lambda: \cf(\delta)=\kappa\}$. If $\cf(\lambda)>\omega$ then $S^\lambda_\kappa$ is a stationary subset of $\lambda$. We shall use the Jerusalem forcing notation, namely $p\leq q$ means that $p$ is weaker than $q$. If $\mathcal{I}$ is an ideal over $\kappa$ then $\mathcal{I}^+ = \mathcal{P}(\kappa)-\mathcal{I}$. We shall always assume that every bounded subset of $\kappa$ belongs to $\mathcal{I}$. The notation ${\rm NS}_\kappa$ refers to the non-stationary ideal over $\kappa$.

Suppose that $\kappa=\cf(\kappa)>\aleph_0$, $S\subseteq\kappa, S$ is stationary and $\Diamond_S$ holds as exemplified by $\langle A_\delta:\delta\in S\rangle$.
Suppose that $M$ is any structure over $\kappa$, and the size of $\mathcal{L}(M)$ is at most $\kappa$.
We would like to say that the diamond sequence predicts elementary submodels of $M$. However, the diamond sequence predicts sets of ordinals, and $M$ contains many objects which are not ordinals.

It is possible to code all the information in $M$ as subsets of $\kappa$.
For this, we fix $|\mathcal{L}(M)|$ disjoint subsets of $\kappa$, each of which is of size $\kappa$, denoted by $\{B_R:R\in\mathcal{L}(M)\}$.
We also fix one-to-one functions from $\kappa^{n(R)}$ into $B_R$ for every $R\in\mathcal{L}(M)$ where $n(R)={\rm arity}(R)$. The union of the range of these functions is a subset of $\kappa$, hence the sequence $\langle A_\delta:\delta\in S\rangle$ predicts it at stationarily many places.

Since the code of each $R^M$ lies in a set disjoint from the other sets and the functions are one-to-one, it is possible to decode the information and recover a submodel of $M$ at each point in which the diamond sequence guesses an initial segment of the above set.
Moreover, the set of ordinals for which such a submodel is elementary will be still a stationary set.
We indicate that the same diamond sequence predicts, in this way, elementary submodels of every structure over $\kappa$.

In the main theorem below, when we say that a certain diamond sequence predicts elementary submodels of some structure $M$, we assume that a coding of the language has been chosen from the outset so that submodels can be deciphered from the sets of ordinals predicted by the diamond sequence.
We will use these models for creating a tiltan sequence with special indestructibility property out of a diamond sequence.

We shall make use of a result of Laver from \cite{MR0472529}, who proved the following. If $\kappa$ is supercompact then there exists a forcing notion which makes $\kappa$ indestructible under any $\kappa$-directed-closed forcing notion.
The forcing of Laver is compatible with having GCH above $\kappa$, and this is important in many applications.
Suppose that $\kappa$ is Laver-indestructible supercompact.
It is possible to force $\mathfrak{s}_\kappa>\kappa^+$ for such a cardinal, and the value of $\mathfrak{s}_\kappa$ can be arbitrarily large. A proof of this fact appears in Claim 2.3 of \cite{MR3201820}, but it was known long before.

\newpage

\section{Tiltan and Galvin's property}

The definition of the tiltan in the previous section asserted the existence of a stationary set of guesses. The following useful lemma says that one inclusion (if guaranteed for every unbounded set) is equivalent to stationarily many inclusions. We spell out the known proof since it is very short.

\begin{lemma}
\label{lemsingle} Suppose that $\kappa=\cf(\kappa)>\aleph_0, S\subseteq\kappa$ and $\langle T_\eta:\eta\in S\rangle$ is a single-guess sequence, i.e. $T_\eta$ is cofinal in $\eta$ for every $\eta\in S$ and for every unbounded $A\subseteq\kappa$ there is an ordinal $\eta\in S$ for which $T_\eta\subseteq A$. \newline
Then the set $\{\eta\in S: T_\eta\subseteq A\}$ is stationary (hence, in particular, $S$ must be a stationary subset of $\kappa$).
\end{lemma}

\par\noindent\emph{Proof}. \newline
Fix an unbounded $A\subseteq\kappa$.
Let $C\subseteq\kappa$ be any club. Suffice it to prove that $T_\eta\subseteq A$ for some $\eta\in S\cap C$.
We do this by defining a set $y\subseteq A$ and a corresponding set $\{\gamma_j:j<\kappa\}\subseteq C$ by induction.
Arriving at stage $i$, we choose $\alpha_i\in A$ such that $\alpha_i>\sup \{\gamma_j:j<i\}$, and then we choose $\gamma_i\in C$ such that $\gamma_i>\alpha_i$. This is possible since both $A$ and $C$ are unbounded in $\kappa$ and $\kappa$ is regular.

Let $y = \{\alpha_i: i<\kappa\}\subseteq A$. By the single-guess assumption we can choose $\eta\in S$ so that $T_\eta\subseteq y\subseteq A$. However, $\eta = \sup(T_\eta)$ and hence $\eta$ is a limit of elements from $C$ since there is such an element betwixt any two members of $y$. As $C$ is closed we see that $\eta\in C$, so we are done.

\hfill \qedref{lemsingle}

Within the proof of the main theorem we shall make use of the strong negation of Galvin's property, as proved consistent in \cite{MR830084}. We quote the first construction of Abraham-Shelah, labeled there as Theorem 1.1. The forcing is done over a model which satisfies GCH.

\begin{theorem}
\label{thmabsh} Abraham-Shelah's theorem. \newline
Suppose that $\kappa=\cf(\kappa)<\kappa^+<\cf(\lambda)\leq\lambda$ and assume GCH. \newline
There is a $\kappa$-complete generic extension which collapses no cardinals, makes $2^{\kappa^+}=\lambda$ and models the following property: \newline
There exists a family of $\lambda$ many clubs of $\kappa^+$ such that the intersection of each $\kappa^+$ of them is of cardinality less than $\kappa$.
\end{theorem}

\hfill \qedref{thmabsh}

In the case $\kappa=\aleph_0$, the above theorem gives finite intersection of any sub-collection of size $\aleph_1$ of the forced family. The following simple claim shows that \emph{finite} cannot be improved to \emph{empty}. A similar statement will be phrased later upon replacing finite by countable, under the assumption of tiltan.

\begin{fact}
\label{clmfinite} For every collection $\{C_\alpha:\alpha<\lambda\}$ of club subsets of $\aleph_1$, where $\lambda=\cf(\lambda)>\aleph_1$, there is a sub-collection of size $\lambda$ with non-empty intersection. Moreover, for every $n\in\omega$ one can find such a sub-collection whose intersection contains at least $n$ elements.
\end{fact}

\par\noindent\emph{Proof}. \newline
Fix any $n\in\omega$, and enumerate $[\aleph_1]^{\geq n}$ by $\{b_\gamma:\gamma\in\omega_1\}$. Let $\{C_\alpha:\alpha<\lambda\}$ be a collection of club subsets of $\aleph_1$. For every $\gamma\in\omega_1$ set $\mathcal{A}_\gamma = \{\alpha\in\lambda: b_\gamma\subseteq C_\alpha\}$. Clearly, $\bigcup_{\gamma\in\omega_1}\mathcal{A}_\gamma = \lambda$.
Since $\lambda=\cf(\lambda)>\aleph_1$, there is an ordinal $\gamma\in\omega_1$ for which $|\mathcal{A}_\gamma|=\lambda$.
It follows that $b_\gamma\subseteq \bigcap\{C_\alpha:\alpha\in \mathcal{A}_\gamma\}$, so we are done.

\hfill \qedref{clmfinite}

The assumption about the regularity of $\lambda$ can be relaxed to $\lambda\geq\cf(\lambda)>\aleph_1$. If $\cf(\lambda)=\aleph_1$ then we can still find a sub-collection of size $\aleph_1$ (or any $\kappa<\lambda$) with the same property.
Assuming now tiltan, we can prove a similar (but stronger) assertion:

\begin{fact}
\label{clmcountable} Assume $\clubsuit$. \newline
For every collection $\{C_\alpha:\alpha\in\omega_2\}$ of clubs of $\aleph_1$ and any ordinal $\tau\in\omega\cdot\omega$ there is a sub-collection of size $\aleph_2$ with intersection of order type at least $\tau$. \newline
More generally, $\kappa=\cf(\kappa)$ and $\clubsuit_{\kappa^+}$ imply that for every family $\{C_\alpha:\alpha\in\kappa^{++}\}$ of clubs of $\kappa^+$ and every ordinal $\tau\in\kappa\cdot\omega$ there is a subfamily of size $\kappa^{++}$ whose intersection has order type $\geq\tau$.
\end{fact}

\par\noindent\emph{Proof}. \newline
Fix a $\clubsuit$-sequence $\langle T_\eta:\eta\in\lim(\omega_1)\rangle$.
Suppose we are given a collection $\{C_\alpha:\alpha\in\omega_2\}$ of clubs of $\aleph_1$ and an ordinal $\tau\in\omega\cdot\omega$.
Choose a natural number $n\in\omega$ so that $\tau\leq\omega\cdot n$.

By induction on $m\leq n$ we choose a set $T_m$ from the tiltan sequence and a collection $\mathcal{F}_m$ of clubs so that $\mathcal{F}_0\subseteq \{C_\alpha:\alpha\in\omega_2\}$, and the following requirements are met:
\begin{enumerate}
\item [$(a)$] $T_m\subseteq \cap\mathcal{F}_m$.
\item [$(b)$] $\ell<m\Rightarrow \mathcal{F}_m\subseteq\mathcal{F}_\ell$.
\item [$(c)$] $|\mathcal{F}_m|=\aleph_2$.
\item [$(d)$] $\ell<m\Rightarrow \sup(T_\ell)<\min(T_m)$.
\end{enumerate}
The choice is possible by a nested application of the pigeon-hole principle.
Indeed, the cardinality of $\mathcal{F}_m$ is $\aleph_2$ at each stage, and every element of $\mathcal{F}_m$ contains some $T_m$ from the tiltan sequence. There are $\aleph_1$-many elements in the tiltan sequence, and after each step $m<n$ we can truncate all the elements of $\mathcal{F}_m$ at $\sup(T_m)+1$ and then choose $T_{m+1}$.

Now let $T = \bigcup_{m\leq n}T_m$, so ${\rm otp}(T)\geq\omega\cdot n\geq\tau$.
For every $C_\alpha\in \mathcal{F}_n$ we can see that $T\subseteq C_\alpha$, so we are done. The additional statement about a larger $\kappa$ can be proved in the same way.

\hfill \qedref{clmcountable}

One may wonder if $\omega\cdot\omega$ is the correct upper bound, and if large intersection can be calibrated at the ordinal level:

\begin{question}
\label{qordinals} Assume $\tau\in\omega_1$ and $\omega\cdot\omega\leq\tau$. \newline
Is it consistent that tiltan holds and there is a collection $\{C_\alpha:\alpha\in\omega_2\}$ of clubs of $\aleph_1$ such that the order type of the intersection of any $\aleph_1$-sub-collection is less than $\tau$?
What about $\aleph_2$-sub-collections?
\end{question}

Actually, we didn't use the closedness of the sets in the above claim, only their unboundedness.
This gives an interesting partition theorem under $\clubsuit$.
Recall that $\binom{\alpha}{\beta} \rightarrow \binom{\gamma}{\delta}^{1,1}_\theta$ means that for every $c:\alpha\times\beta \rightarrow\theta$ one can find $A\subseteq\alpha, {\rm otp}(A)=\gamma,$ and $B\subseteq\beta, {\rm otp}(B)=\delta,$ such that $c\upharpoonright(A\times B)$ is constant.

\begin{corollary}
\label{corpartition} Assume $\clubsuit$.
\begin{enumerate}
\item [$(\aleph)$] For every family $\{A_\alpha:\alpha\in\omega_2\} \subseteq [\aleph_1]^{\aleph_1}$ and any $\tau\in\omega\cdot\omega$ there is $B\subseteq\omega_2, |B|=\aleph_2$ such that ${\rm otp}(\bigcap\{A_\alpha:\alpha\in B\})\geq\tau$.
\item [$(\beth)$] $\binom{\omega_2}{\omega_1} \rightarrow \binom{\omega_2}{\tau}^{1,1}_{\aleph_0}$ for every $\tau\in\omega\cdot\omega$.
\item [$(\gimel)$] The relation $\binom{\omega_2}{\omega_1} \rightarrow \binom{\omega_2}{\tau}^{1,1}_{\aleph_0}$ for every $\tau\in\omega\cdot\omega$ is consistent with both $\binom{\omega_2}{\omega} \rightarrow \binom{\omega_2}{\omega}^{1,1}_2$ and $\binom{\omega_2}{\omega} \nrightarrow \binom{\omega_2}{\omega}^{1,1}_2$.
\end{enumerate}
\end{corollary}

\par\noindent\emph{Proof}. \newline
Part $(\aleph)$ follows from the proof of Theorem \ref{clmcountable} as indicated above. For part $(\beth)$, let $c:\omega_2\times\omega_1\rightarrow \aleph_0$ be a coloring. For every $\alpha\in\omega_2, n\in\omega$ let $A_{\alpha n} = \{\beta\in\omega_1: c(\alpha,\beta)=n\}$.
For each $\alpha\in\omega_2$ there is $n(\alpha)\in\omega$ so that $A_{\alpha n(\alpha)}$ is uncountable.
We may assume without loss of generality that $n(\alpha)=n$ for some fixed $n$, for every $\alpha\in\omega_2$.
By part $(\aleph)$ there is $B\subseteq\omega_2, |B|=\aleph_2$ such that $y = \bigcap\{A_{\alpha n}:\alpha\in B\}$ satisfies ${\rm otp}(y)\geq\tau$.
Since $c\upharpoonright(B\times y)$ is constantly $n$, we are done.

If $\mathfrak{c}=\omega_2$ then $\binom{\omega_2}{\omega} \rightarrow \binom{\omega_2}{\omega}^{1,1}_2$ iff $\mathfrak{r}=\aleph_1$ (see Claim 1.4 of \cite{MR3201820} for one direction and Main Claim 1.1 of \cite{1047} for the other direction).
Now in Baumgartner's model of side-by-side Sacks forcing we have that $\mathfrak{r}=\aleph_1$ and $\clubsuit$ holds (see \cite{MR2279653}). On the other hand, Brendle \cite{MR2279653} proved the consistency of $\clubsuit$ with $\mathfrak{c}={\rm cov}(\mathcal{B})=\omega_2$.
Here, $\mathcal{B}$ stands for Baire and ${\rm cov}(\mathcal{B})$ is therefore the covering number of the ideal of meager sets.
In this model, $\mathfrak{r}=\omega_2$ as well since ${\rm cov}(\mathcal{B})\leq \mathfrak{r}$. This proves both directions of part $(\gimel)$, so the proof is accomplished.

\hfill \qedref{corpartition}

The moral of the above corollary is that tiltan is inconsistent with the strong negation of Galvin's property.
Actually, the corollary proves something stronger.
It shows that under tiltan one can find sub-collections of size $\aleph_2$ with countable intersection.
Likewise, this can be done for any family of unbounded sets in $\omega_1$, even if these sets are not clubs.
In particular, tiltan fails in the models of \cite{MR830084}.
Nevertheless, Galvin's property does not follow from tiltan.

\begin{theorem}
\label{thmneg} Tiltan is consistent with the negation of Galvin's property.
\end{theorem}

\par\noindent\emph{Proof}. \newline
We force over $L$ as our ground model, so $\Diamond_S$ holds at every stationary $S\subseteq\kappa=\cf(\kappa)>\aleph_0$.
We shall use it in order to construct a $\clubsuit_S$-sequence for some $S\subseteq S^{\aleph_2}_{\aleph_0}$.
The special property of this sequence will be indestructibility under $\aleph_1$-complete forcing notions.
The definition of the sequence and the indestructibility proof are elaborated in \cite{MR1623206}, I, \S 7.
For completeness, we describe this construction explicitly.

We expand the language of set theory by the relation symbols $P$ and $R$, both of them are two-placed relations. We shall concentrate on structures on $\aleph_2$ in which $P$ is interpreted as a partial order, $P(\beta,\gamma)$ implies that for every $\delta$, $R(\beta,\delta)\rightarrow R(\gamma,\delta)$ is true, and for all $\alpha$ and $\beta$ there are $\gamma$ and $\delta$ such that $\alpha<\gamma, \beta<\delta$ and $R(\gamma,\delta)$.

Let $(A_\eta:\eta\in S^{\aleph_2}_{\aleph_0})$ be a diamond sequence, and let $M = (\aleph_2,P,R)$ be a structure for which the relations $P$ and $R$ satisfy the above requirements.
Let $S\subseteq S^{\aleph_2}_{\aleph_0}$ be stationary, and assume that $\eta\in S$ implies that $A_\eta$ guesses $M\upharpoonright\eta$ in the sense described in the introduction. Let $M_\eta = (\eta,P_\eta,R_\eta)$ be the elementary submodel of $M$ deciphered from $A_\eta$.

Our tiltan sequence $(T_\eta:\eta\in S)$ will be constructed as follows.
Firstly, we choose a cofinal sequence of ordinals $\langle\beta^\eta_n:n\in \omega\rangle$ which converges to $\eta$, for every $\eta\in S$. If $A_\eta$ encodes $M_\eta=(\eta,P_\eta,R_\eta)$, then we choose, by induction on $n\in\omega$, a pair of ordinals $(\gamma^\eta_n,\delta^\eta_n)$ such that:
\begin{enumerate}
\item [$(\aleph_n)$] $P_\eta(\gamma^\eta_n,\gamma^\eta_{n+1})$.
\item [$(\beth_n)$] $\delta^\eta_n>\beta^\eta_n$.
\item [$(\gimel_n)$] $R_\eta(\gamma^\eta_n,\delta^\eta_n)$.
\end{enumerate}
We define $T_\eta = \{\delta^\eta_n:n\in\omega\}$ for every $\eta\in S$. Notice that $T_\eta\subseteq\eta$ and by $(\beth_n)$ for every $n\in\omega$ it is cofinal in $\eta$. Our tiltan sequence is $(T_\eta:\eta\in S)$.

For showing its predictive ability, assume that $u$ is an unbounded subset of $\aleph_2$. Let $M$ be the structure whose universe is $\omega_2$, $P^M$ is the order of the ordinals and $R^M = \{(a,b):b\in u\}$. Let $S_0\subseteq S$ be stationary such that $M_\eta = M\cap\eta\prec M$ for every $\eta\in S_0$. It follows from the definition of $R^M$ and elementarity that $T_\eta\subseteq u$ for every $\eta\in S_0$.
We indicate that this specific tiltan construction is shown in \cite{MR1623206} to be indestructible under $\aleph_1$-complete forcing notions, hence we can proceed as follows.

Choose $\lambda\geq\cf(\lambda)>\aleph_2$, and force with the forcing from Theorem \ref{thmabsh} in order to add $\lambda$ many club subsets of $\aleph_2$, the intersection of any sub-collection of which of size $\aleph_2$ is bounded in $\aleph_2$. This forcing is $\aleph_1$-complete and hence $\langle T_\eta:\eta\in S\rangle$ remains a tiltan sequence in the generic extension.
By properness, $S$ will be stationary in the generic extension. It follows that `the champion of the day remained stationary' (see \cite{Ivanhoe}, p. 82).

Call the above forcing $\mathbb{P}$.
We make the comment that it preserves stationary subsets of $\aleph_2$, though this is not need for our argument.
All we need is that $\mathbb{P}$ preserves the stationary sets on which $\langle T_\eta:\eta\in S\rangle$ guesses.
Moreover, Lemma \ref{lemsingle} shows that one guess implies stationarily many guesses, hence if we prove the single guess property then we are done (and this will show, in particular, that the domain of the tiltan is still stationary).

Fix a generic subset $G$ of $\mathbb{P}$.
Let $\mathbb{Q}$ be the usual L\'evy collapse of $\aleph_1$ to a countable ordinal, as defined in $V[G]$.
Let $H\subseteq\mathbb{Q}$ be generic over $V[G]$.
The key-point of this stage is that if $\name{B}$ is a $\mathbb{Q}$-name in $V[G]$ of a set whose $H$-interpretation is an unbounded subset of $\aleph_1^{V[G][H]}$ then there are $A\in V[G], A$ unbounded in $\aleph_2^{V[G]}$ and $q\in\mathbb{Q}$ so that $q \Vdash A\subseteq\name{B}$.
This point will ensure both $\clubsuit_{\aleph_1}$ and the negation of Galvin's property at $\aleph_1$ in $V[G][H]$.

Let us prove the key-point.
Suppose that $\name{B}$ is as above and choose any condition $p$ which forces $\name{B}$ to be unbounded in $\aleph_1^{V[G][H]}$.
Recall that $\aleph_2^{V[G]}$ is collapsed to be $\aleph_1^{V[G][H]}$, so working in $V[G]$ we choose for every $i\in\omega_2$ a condition $q_i\in\mathbb{Q}$ and an ordinal $\tau_i\in\omega_2$ such that:
\begin{enumerate}
\item [$(a)$] $p\leq q_i$.
\item [$(b)$] $i<\tau_i$.
\item [$(c)$] $q_i\Vdash\tau_i\in\name{B}$.
\end{enumerate}
Since $V[G]\models|\mathbb{Q}|=\aleph_1$, we can choose $E\subseteq\omega_2, |E|=\aleph_2$, and a fixed condition $q\in\mathbb{Q}$ such that $i\in E \Rightarrow q_i = q$.
Set $A = \{\tau_i:i\in E\}$ and notice that $q\Vdash A\subseteq\name{B}$. Since $A$ is unbounded in $\omega_2$ and $A\in V[G]$, we are done proving the key-point.

Back to the tiltan sequence $\langle T_\eta:\eta\in S\rangle$ in $V[G]$, we claim that it is still a tiltan sequence in $V[G][H]$, but now for subsets of $\aleph_1$.
Indeed, suppose that $B\in V[G][H]$ is an unbounded subset of $\aleph_1$, and let $\name{B}$ be a name of $B$ in $V[G]$.
By the key-point, choose $A\in V[G]$ so that $A$ is unbounded in $\aleph_2^{V[G]}$ and $q \Vdash_{\mathbb{Q}} \check{A}\subseteq\name{B}$ for some $q\in\mathbb{Q}$.
In $V[G]$ we have an ordinal $\eta\in S$ so that $T_\eta\subseteq A$, so $V[G][H]\models T_\eta\subseteq A\subseteq B$ and hence $\langle T_\eta:\eta\in S\rangle$ is a tiltan sequence by virtue of Lemma \ref{lemsingle}.

For the negation of Galvin's property, fix in $V[G]$ a collection $\{C_\alpha:\alpha<\lambda\}$ of club subsets of $\aleph_2$ which exemplifies this statement.
Namely, any sub-collection of it of size $\aleph_2$ has bounded intersection in $\aleph_2$.
In $V[G][H]$ this is a collection of at least $\aleph_2$-many club subsets of $\aleph_1$ (recall that $\lambda^{V[G]}>\aleph_2^{V[G]}$ and $\lambda$ remains a cardinal after forcing with $\mathbb{Q}$).
If we had a sub-collection of size $\aleph_1$ whose intersection contains an unbounded set $B$ (hence a club) in $\aleph_1$, then this intersection would have a name $\name{B}\in V[G]$.
In $V[G]$ we could find a condition $q$ and some $A\subseteq\aleph_2^{V[G]}$ unbounded, such that $q \Vdash_{\mathbb{Q}}A\subseteq\name{B}$.
But then, the corresponding sub-collection in $V[G]$ will be of size $\aleph_2$ and $A$ will be a subset of its intersection, which is impossible.

\hfill \qedref{thmneg}

Let us take a broader look at $\clubsuit$ and $\Diamond$.
Assume that $\varphi$ is a statement in the language of set theory.
Suppose, further, that the diamond implies $\varphi$, and we wish to force $\clubsuit + \neg\varphi$.
The method of \cite{MR1623206} makes this possible, provided that $\neg\varphi$ can be forced over $\aleph_2$ (or some larger cardinal) by an $\aleph_1$-complete forcing notion, and be preserved by the pertinent collapse.

One may wonder if this is the only way to force tiltan along with the negation of statements which follow from diamond.
To be more concrete, consider the statement $\mathfrak{s}=\omega_1$ which follows from the diamond (or even CH).
The possible consistency of tiltan with its negation, namely $\mathfrak{s}\geq\omega_2$, is an open problem. Apparently, this statement is not amenable to the above properties of $\varphi$.

ZFC implies that $\mathfrak{s}_{\aleph_1}=\aleph_0$, since $\mathfrak{s}_\kappa\geq\kappa$ iff $\kappa$ is strongly inaccessible.
Hence the attempt to force $\mathfrak{s}_\kappa>\kappa^+$ and then to collapse cardinals in order to get $\mathfrak{s}>\omega_1$ with tiltan will fail.
Still, one may wish to force $\mathfrak{s}_\kappa>\kappa^+$ at a very large cardinal $\kappa$, and this requires that $\kappa$ be at least weakly compact (probably much more).
However, the collapse of $\kappa$ to $\aleph_0$ will of course destroy this property.
The pattern of proof used in this section seems to be effective only when forcing the desired property over some $\aleph_n$.
Based on these considerations we conjecture that tiltan implies $\mathfrak{s}=\aleph_1$.

In the second section of \cite{MR830084} the authors force a strong negation of Galvin's property (namely, with finite intersections), which is indestructible under any extension which does not collapse $\aleph_1$ and $\aleph_2$.
It follows from Theorem \ref{clmcountable} that any generic extension of the universe obtained by \cite{MR830084} (in the second section) in which $\clubsuit$ is forced, collapses $\aleph_1$ or $\aleph_2$.
Of course, if one begins with a different ground model then tiltan can be forced without any collapse.

This gives an insight to both the possible ways to force tiltan and the possible statements which are consistent with it. In order to demonstrate the above idea we introduce a negative answer to Question 9.0.26 from \cite{primavesi}.
Primavesi asked whether $\clubsuit+\neg{\rm CH}$ implies $2^\omega=2^{\omega_1}$. Using the methods of this section, we can prove the following.

\begin{theorem}
\label{tilprimavesi} Tiltan and weak diamond. \newline
It is consistent that tiltan holds, the continuum hypothesis fails and $2^\omega<2^{\omega_1}$. Moreover, the values of $2^\omega$ and $2^{\omega_1}$ can be arbitrarily large, as well as the discrepancy between them.
\end{theorem}

\par\noindent\emph{Proof}. \newline
Begin with a stationary $S\subseteq S^{\aleph_2}_{\aleph_0}$ and a $\clubsuit_S$-sequence which is indestructible upon $\aleph_1$-complete forcing notions.
Choose $\lambda\geq\cf(\lambda)>\aleph_2$ and $\kappa\geq\cf(\kappa)>\aleph_1$ so that $\kappa<\lambda$, and force $2^{\omega_1}=\kappa<\lambda=2^{\omega_2}$ using any $\aleph_1$-complete forcing notion (e.g., by adding Cohen subsets to $\omega_1$ and $\omega_2$).
Finally, collapse $\aleph_1$ to $\aleph_0$.
The resulting model satisfies the statement of the theorem.

\hfill \qedref{tilprimavesi}

Let us conclude with one more statement which seems to be connected with Galvin's property.
Recall that an ideal $\mathcal{I}$ over $\aleph_1$ is saturated iff in every collection $\{S_\alpha:\alpha\in\omega_2\}\subseteq\mathcal{I}^+$ there are two distinct sets $S_\alpha$ and $S_\beta$ so that $S_\alpha\cap S_\beta\in \mathcal{I}^+$. Let us focus on $\mathcal{I} = {\rm NS}_{\omega_1}$.

The general impression is that Galvin's property implies that ${\rm NS}_{\omega_1}$ is not saturated. Intuitively, if the club subsets are not too scattered (in the sense given by Galvin's property) then one can find many stationary sets which are quite separated. A basic example is under the diamond, when Galvin's property holds and there are $\aleph_2$-many stationary sets which are almost disjoint. This means, of course, that the non-stationary ideal is not saturated. This is the case also under tiltan:

\begin{claim}
\label{clmad} If $\clubsuit$ holds then there is an almost disjoint family of stationary subsets of $\aleph_1$ of size $\aleph_2$, hence the negation of Galvin's property is consistent with ${\rm NS}_{\omega_1}$ being non-saturated.
\end{claim}

\par\noindent\emph{Proof}. \newline
Let $\langle T_\eta:\eta\in\omega_1\rangle$ be a tiltan sequence.
Let $\{A_\beta:\beta\in\omega_2\}$ be an almost disjoint family of unbounded subsets of $\aleph_1$.
For each $\beta\in\omega_2$ let $S_\beta$ be the set $\{\eta\in\omega_1: T_\eta\subseteq A_\beta\}$. So every $S_\beta$ is a stationary subset of $\aleph_1$, and the family $\{S_\beta:\beta\in\omega_2\}$ must be almost disjoint.
Hence tiltan implies that ${\rm NS}_{\omega_1}$ is not saturated.
In particular, if we force tiltan with the negation of Galvin's property then we prove the second part of the claim.

\hfill \qedref{clmad}

The above claim invites a natural question:

\begin{question}
\label{qmm} Martin's Maximum and Galvin's property. \newline
Does the negation (or strong negation) of Galvin's property follow from Martin's Maximum?
\end{question}

\newpage

\section{Superclub and Galvin's property}

In the previous section we have seen that tiltan is consistent with Galvin's property.
This requires a violation of CH, since under the continuum hypothesis Galvin's property holds.
It is clear that Galvin's property follows from diamond, since diamond implies CH. The purpose of this section is to show that the crucial point is not CH but rather the prediction element of the diamond.
We do this by analyzing the situation under an intermediate prediction principle between diamond and tiltan.

\begin{definition}
\label{defsuperclub} Superclub.
\begin{enumerate}
\item [$(\aleph)$] A sequence of sets $\langle A_\delta:\delta\in\lim(\omega_1)\rangle$ is a superclub sequence iff each $A_\delta$ is a cofinal subset of $\delta$ and for every unbounded $x\subseteq\omega_1$ there is an unbounded $y\subseteq x$ such that the set $\{\delta\in\omega_1:y\cap\delta=A_\delta\}$ is a stationary subset of $\omega_1$.
\item [$(\beth)$] Superclub (at $\aleph_1$) means that there exists a superclub sequence (of legnth $\omega_1$).
\end{enumerate}
\end{definition}

Superclub has been defined in \cite{primavesi}, where the author shows that it implies the existence of a Suslin tree.
It is clear from the definition that diamond implies superclub and superclub implies tiltan.
Both implications are irreversible.
The consistency of superclub with the negation of diamond (and even the negation of the continuum hypothesis) is proved in \cite{MR3687435}.
The consistency of tiltan with the negation of superclub is easier, and will be demonstrated also in the main theorem of this section.
Namely, we shall prove that superclub implies Galvin's property.
Ahead of the proof we need the following statement from \cite{primavesi}.

\begin{lemma}
\label{lemclubclub} Assume superclub. \newline
Then there exists a sequence $\langle B_\delta:\delta\in\lim(\omega_1)\rangle$ such that every $B_\delta$ is a club susbet of $\delta$ and for every club $C$ of $\omega_1$ there is a club $D\subseteq C$ such that the set $\{\delta\in\omega_1:D\cap\delta=B_\delta\}$ is stationary.
\end{lemma}

\par\noindent\emph{Proof}. \newline
Let $\langle A_\delta:\delta\in\lim(\omega_1)\rangle$ be a superclub sequence.
For every $\delta\in\lim(\omega_1)$ let $B_\delta$ be the closure of $A_\delta$ in the order topology.
We claim that the sequence $\langle B_\delta:\delta\in\lim(\omega_1)\rangle$ is as required.
For this, fix a club subset $C$ of $\omega_1$.
Let $y$ be an unbounded subset of $C$ such that $\{\delta\in\omega_1:y\cap\delta=A_\delta\}$ is stationary in $\omega_1$.

Set $D=c\ell(y)$, so $D$ is a club of $\omega_1$.
If $\gamma\in\{\delta\in\omega_1:y\cap\delta=A_\delta\}$ then $A_\gamma\subseteq y$ and hence $B_\gamma\subseteq D$.
Inasmuch as $y\subseteq C$ and $C$ is closed we see that $D\subseteq C$.
Since $\{\delta\in\omega_1:y\cap\delta=A_\delta\}$ is stationary we conclude that $\{\gamma\in\omega_1:D\cap\gamma=B_\gamma\}$ is stationary.

\hfill \qedref{lemclubclub}

We shall call a sequence of the form $\langle B_\delta:\delta\in\lim(\omega_1)\rangle$ as guaranteed by the lemma \emph{a closed superclub sequence}.

\begin{theorem}
\label{thmgp} Superclub implies Galvin's property.
\end{theorem}

\par\noindent\emph{Proof}. \newline
Fix a closed superclub sequence $\langle B_\delta:\delta\in\lim(\omega_1)\rangle$.
Let $\{C_\alpha:\alpha\in\omega_2\}$ be any collection of club subsets of $\omega_1$.
For every $\alpha\in\omega_2$ we choose a club $D_\alpha$ of $\omega_1$ so that $D_\alpha\subseteq C_\alpha$ and the set $\{\delta\in\omega_1: D_\alpha\cap\delta=B_\delta\}$ is stationary.
For every $\alpha\in\omega_2, \delta\in\omega_1$ we define the set $H_{\alpha\delta}$ as follows.
If $D_\alpha\cap\delta\neq B_\delta$ then $H_{\alpha\delta}=\varnothing$.
If $D_\alpha\cap\delta=B_\delta$ then define:
$$
H_{\alpha\delta}=\{\beta\in\omega_2: D_\beta\cap\delta=B_\delta\}.
$$
Notice that for every $\alpha\in\omega_2$ there is an unbounded (and even stationary set) of $\delta\in\omega_1$ for which $H_{\alpha\delta}$ is not empty.
Notice, further, that if $\gamma<\delta$ and both $H_{\alpha\gamma}$ and $H_{\alpha\delta}$ are not empty then $H_{\alpha\delta}\subseteq H_{\alpha\gamma}$.

For this, assume that $\beta\in H_{\alpha\delta}$.
By definition, $D_\beta\cap\delta=D_\alpha\cap\delta=B_\delta$.
Hence $D_\beta\cap\gamma = (D_\beta\cap\delta)\cap\gamma = (D_\alpha\cap\delta)\cap\gamma = D_\alpha\cap\gamma$.
But $D_\alpha\cap\gamma=B_\gamma$ since $H_{\alpha\gamma}$ is not empty.
It follows that $D_\beta\cap\gamma=D_\alpha\cap\gamma=B_\gamma$ so $\beta\in H_{\alpha\gamma}$ as required.

As in the original proof of Galvin under the continuum hypothesis, we can argue that for some $\alpha\in\omega_2$ it it true that $|H_{\alpha\delta}|=\aleph_2$ for an unbounded set of $\delta\in\omega_1$.
If not, then for every $\alpha\in\omega_2$ there is a first ordinal $\delta_\alpha\in\omega_1$ such that $|H_{\alpha\delta_\alpha}|<\aleph_2$.
By Fodor's lemma, there are a stationary subset $S\subseteq\omega_2$ and a fixed $\delta\in\omega_1$ such that $\alpha\in S\Rightarrow \delta_\alpha=\delta$.
Observe that if $H_{\alpha\delta}\neq\varnothing$ then it does not depend on $\alpha$ since it is determined solely by $B_\delta$.
Since $H_{\alpha\delta}\neq\varnothing$ for every $\alpha\in S$ we see that $H_{\alpha\delta}=T$ for some fixed $T$ of size less than $\aleph_2$ by our assumption.
It follows that $\bigcup_{\alpha\in S}H_{\alpha\delta}=T$.
However, $\alpha\in S\Rightarrow\alpha\in T$, so $|T|\geq|S|=\aleph_2$, a contradiction.

Fix, therefore, an ordinal $\alpha\in\omega_2$ for which $|H_{\alpha\delta}|=\aleph_2$ whenever $H_{\alpha\delta}$ is not empty.
Enumerate the non-empty $H_{\alpha\delta}$s by $\{H_{\alpha\delta_\xi}:\xi\in\omega_1\}$.
By induction on $\xi\in\omega_1$ choose $\beta_\xi\in H_{\alpha\delta_\xi}$ such that $\alpha_\xi\notin\{\alpha_\eta:\eta<\xi\}$.
It is easily verified that the intersection of  $\{D_{\beta_\xi}:\xi\in\omega_1\}$ is a club subset of $\omega_1$.
This is true also for $\{C_{\beta_\xi}:\xi\in\omega_1\}$, so we are done.

\hfill \qedref{thmgp}

In the light of the previous section, the above theorem shows that tiltan is consistent with the negation of superclub.
Overall, the strong negation of Galvin's property is consistent but tiltan would fail in such models. Tiltan is still consistent with the weak negation of Galvin's property, and superclub is stronger in the sense that it implies Galvin's property.
We point to the reason that superclub cannot be replaced by tiltan in the above proof.
If $\langle B_\delta:\delta\in\lim(\omega_1)\rangle$ is just a tiltan sequence then the fact that $\beta\in H_{\alpha\delta}$ does not guarantee $\beta\in H_{\alpha\gamma}$ when $\gamma<\delta$ and $H_{\alpha\gamma}\neq\varnothing$.

\newpage

\section{Tiltan and cardinal characteristics}

A cardinal characteristic is defined, usually, as the minimal size of a family with some property. Cardinal characteristics assume a value in the interval $[\omega_1,2^\omega]$, and if not trivial they admit a variety of possible values.
These values are limited, in most cases, by two types of constraints.
Firstly, the possible cofinality of the given characteristic. Secondly, its relation with other characteristics.
These relations are usually amenable to consistency results.

One can draw a diagram in which there are several robust ZFC relations and some freedom. In such diagrams, the smallest nicely defined characteristic is $\mathfrak{m}$, the minimal $\kappa$ for which ${\rm MA}_\kappa$ fails. It is known that $\mathfrak{m}\leq\mathfrak{p}$ and consistently $\mathfrak{m}<\mathfrak{p}$, so $\mathfrak{m}$ is small indeed.
From the other side, the independence number $\mathfrak{i}$ is one of the largest characteristics. In fact, there is no other nicely defined characteristic which is, provably in ZFC, at least $\mathfrak{i}$ (and most other characteristics are always below $\mathfrak{i}$).

We shall define, below, a cardinal characteristic dubbed as $\mathfrak{gp}$.
It will be shown that $\mathfrak{gp}$ is totally free when considering the relationship with other classical characteristics.
In particular, $\mathfrak{gp}<\mathfrak{m}$ is consistent, as well as $\mathfrak{gp}>\mathfrak{i}$.

\begin{definition}
\label{defgp} The cardinal characteristic $\mathfrak{gp}$. \newline
We define $\mathfrak{gp}$ as the minimal $\kappa$ such that every family $\{C_\alpha:\alpha\in\kappa^+\}$ of club subsets of $\aleph_1$ contains a subfamily $\{C_{\alpha_\beta}:\beta\in\omega_1\}$ whose intersection is closed and unbounded in $\aleph_1$.
\end{definition}

The name $\mathfrak{gp}$ comes from the fact that the collection mentioned in the definition satisfies Galvin's property.
We commence with the following simple statements:

\begin{proposition}
\label{propbasics} Basic facts.
\begin{enumerate}
\item [$(\aleph)$] $\mathfrak{gp}$ is well-defined.
\item [$(\beth)$] $\aleph_1\leq\mathfrak{gp}\leq 2^{\aleph_0}$.
\end{enumerate}
\end{proposition}

\par\noindent\emph{Proof}. \newline
Let $\kappa=2^\omega$.
Assume that $\{C_\alpha:\alpha\in\kappa^+\}$ is a collection of club subsets of $\aleph_1$, we may allow repetitions (so actually, this is a sequence).
If $2^{\omega_1}=\kappa$ as well then there is a club $C$ which appears $\kappa^+$-many times in the above collection, so clearly this family satisfies Galvin's property.
If $2^\omega=\kappa<2^{\omega_1}$ then $\{C_\alpha:\alpha\in\kappa^+\}$ satisfies Galvin's property by Theorem 2.1 of \cite{MR3604115}. In both cases, every family of clubs of size $\kappa^+$ satisfies Galvin's property and hence $\mathfrak{gp}$ is well-defined.

For the second part of the proposition notice that the above argument gives $\mathfrak{gp}\leq 2^\omega$.
Now for every $\alpha\in\omega_1$ let $C_\alpha = \{\beta\in\omega_1:\beta>\alpha\}$ and let $\mathcal{F} = \{C_\alpha:\alpha\in\omega_1\}$.
So $\mathcal{F}$ is a collection of $\aleph_1$ clubs of $\aleph_1$, and any uncountable subfamily of $\mathcal{F}$ has empty intersection.
This means that $\mathfrak{gp}>\aleph_0$, so we are done.

\hfill \qedref{propbasics}

Ahead of the theorem below we mention that the negation of Galvin's property is forced in Theorem 2.2 of \cite{MR830084} with an additional property. For $\lambda\geq\aleph_2$, the authors force a collection $\{C_\alpha:\alpha<\lambda\}$ of clubs of $\aleph_1$, any sub-collection of which of size $\aleph_1$ has finite intersection and this family keeps its property in any generic extension which preserves $\aleph_1$.
There is something misleading in the last statement, since if $\lambda=\aleph_2$ and $\aleph_2$ is collapsed then the family keeps the property of the negation of Galvin's property but its size becomes $\aleph_1$ so it is meaningless.
Anyway, in any generic extension which collapses no cardinals, there is still a family $\{C_\alpha:\alpha<\lambda\}$ of clubs of $\aleph_1$ which negates Galvin's property.

\begin{theorem}
\label{thmfinal} Assume GCH and suppose that $\theta=\cf(\theta)>\aleph_1$.
\begin{enumerate}
\item [$(a)$] For every $\lambda>\cf(\lambda)=\theta$ there is a cardinal-preserving generic extension in which $\mathfrak{gp}=\lambda$.
\item [$(b)$] It is consistent that $\mathfrak{gp}<\mathfrak{m}$, and the gap can be arbitrarily large.
\item [$(c)$] It is consistent that $\mathfrak{gp}>\mathfrak{i}$, and the gap can be arbitrarily large.
\end{enumerate}
\end{theorem}

\par\noindent\emph{Proof}. \newline
Begin with a model of GCH, and force with the forcing given by Theorem \ref{thmabsh} which makes $2^{\aleph_1}=\lambda$ and produces a collection $\{C_\alpha:\alpha\in\lambda\}$ of clubs of $\aleph_1$ without Galvin's property.
Notice that $2^\omega=\lambda$ as well, by Theorem 2.1 of \cite{MR3604115}. It follows that $\lambda\leq\mathfrak{gp}\leq 2^\omega=\lambda$, so part $(a)$ is proved.
We observe that this gives the consistency of $\mathfrak{g}<\mathfrak{gp}$ or $\mathfrak{b}<\mathfrak{gp}$, or any other characteristic which is always regular.

Part $(b)$ follows from Theorem 2.5 of \cite{MR3604115}, but we unfold the argument.
If $\mathcal{F} = \{C_\alpha:\alpha<\lambda\}$ satisfies Galvin's property, $\mathbb{P}$ is a $ccc$ forcing notion and $G\subseteq\mathbb{P}$ is generic then $\mathcal{F}$ still satisfies Galvin's property in $V[G]$.
The reason is that every club of $\aleph_1$ in the generic extension contains an old club.
In particular, if the ground model satisfies $2^\omega=\omega_1$ (so $\mathfrak{gp}=\aleph_1$) and $\mathbb{P}$ forces ${\rm MA} + 2^\omega=\lambda$ in the usual way, then $V[G]\models \mathfrak{gp}=\aleph_1<\lambda = \mathfrak{m}$.

Finally, we choose any prescribed $\lambda\geq\aleph_2$ and we force with the second forcing described in \cite{MR830084}, in order to get $\mathfrak{gp}\geq\lambda$ by a forcing notion $\mathbb{P}$. Let $G\subseteq\mathbb{P}$ be generic over $V$.
Working in $V[G]$ we can choose any forcing notion $\mathbb{Q}$ which preserves cardinals (or even just $\aleph_1^{V[G]}$ and $\aleph_2^{V[G]}$), and we can take a subset $H\subseteq\mathbb{Q}$ which is $V[G]$-generic.
From the special property of the set which witnesses the negation of Galvin's property added by $\mathbb{P}$ we infer that $V[G][H]\models \mathfrak{gp}\geq\lambda$.
This holds, in particular, whenever $\mathbb{Q}$ is $ccc$ in $V[G]$.

One can force, in this way, $\mathfrak{i}=\omega_1$ without changing the value of $2^\omega$ in $V[G]$. Hence $\mathfrak{i}<\mathfrak{gp}$ holds in $V[G][H]$. The same holds for other statements which can be forced by a $ccc$ forcing notion. If such a forcing diminishes specific characteristics then the large value of $\mathfrak{gp}$ is preserved in the generic extension.

\hfill \qedref{thmfinal}

One may wonder about small cofinality:

\begin{question}
\label{qcof} Is it provable that $\cf(\mathfrak{gp})>\omega$?
\end{question}

The main theorem of the first section can be described now as the consistency of tiltan with $\mathfrak{gp}=2^\omega$.
It is known that tiltan implies $\mathfrak{p}=\omega_1$, but even for a small characteristic like $\mathfrak{h}$ it is unknown whether tiltan is consistent with $\mathfrak{h}>\omega_1$. Large characteristics can be forced to be above $\aleph_1$ along with tiltan, but it seems that the reason for the consistency of tiltan with a large value of $\mathfrak{gp}$ is different.

Almost all cardinal characteristics are connected with a property of subsets of $\omega$, while $\mathfrak{gp}$ is basically connected with subsets of $\aleph_1$. This is apparently the reason for the extreme freedom of $\mathfrak{gp}$ with respect to other characteristics, as well as the result concerning tiltan.

The possibility of $\clubsuit + \mathfrak{s}>\omega_1$ is still open. The difficulty was mentioned in the previous section. Specifically, the method of working at the $\aleph_2$-level and then collapsing $\aleph_1$ is not helpful with respect to the splitting number.
It seems, however, that this idea is applicable for $\mathfrak{s}_\kappa$ at large cardinals, and this is the main theorem of this section. There are several theorems about cardinal characteristics which hold (in ZFC) at uncountable $\kappa$ but fail when $\kappa=\aleph_0$.
A salient example is the statement $\mathfrak{s}_\kappa\leq\mathfrak{b}_\kappa$ (proved in \cite{MR3615051} for every $\kappa>\omega$).
In the classical case of $\kappa=\aleph_0$ it is consistent that $\mathfrak{b}<\mathfrak{s}$.

\begin{theorem}
\label{thmsupercompact} Let $\kappa$ be a supercompact cardinal. \newline
One can force $\clubsuit_{\kappa^+}$ along with $\mathfrak{s}_\kappa>\kappa^+$, while preserving the supercompactness of $\kappa$.
\end{theorem}

\par\noindent\emph{Proof}. \newline
We begin with the preparatory forcing of Laver which makes $\kappa$ indestructible under $\kappa$-directed-closed forcing notions, and we also force $2^\kappa=\kappa^+$.
As a second step we create a tiltan sequence of the form $\langle T_\eta:\eta\in S^{\kappa^+}_{\omega}\rangle$ such that ${\rm otp}(T_\eta)=\omega$ for every $\eta\in S^{\kappa^+}_{\omega}$ and the tiltan is indestructible by any $\aleph_1$-complete forcing notion.
For this construction it is enough to assume that $2^\kappa=\kappa^+$ in the ground model, so $\Diamond_{S^{\kappa^+}_{\omega}}$ holds as proved in \cite{MR2596054}. The construction proceeds now exactly as in the first section upon replacing $\aleph_2$ there by $\kappa^+$ here.
A similar argument to the $\aleph_2$-case (which appears in \cite{MR1623206}) shows that this tiltan sequence is indestructible with respect to $\aleph_1$-complete forcing notions, and we spell out the proof.

Let $\mathbb{P}$ be an $\aleph_1$-complete forcing notion.
Let $\name{a}$ be a $\mathbb{P}$-name of an unbounded subset of $\kappa^+$, and let $p$ be any condition in $\mathbb{P}$ which forces it. Fix a sufficiently large regular $\chi>\kappa^+$ and choose an elementary submodel $N\prec\mathcal{H}(\chi)$ such that $p, \mathbb{P}, \name{a}\in N$ and $|N|=\kappa^+$.

Enumerate the elements of $\mathbb{P}\cap N$ above $p$ by $\{p_i: i\in\kappa^+\}$. Using this enumeration we define another structure, $M$, as follows. We add the relation symbols $P$ and $R$ to $\mathcal{L}(M)$. The universe of $M$ will be $\kappa^+$. We shall interpret $P$ by saying that $P(i,j)$ iff $p_i\leq_{\mathbb{P}} p_j$ and $R(i,\delta)$ iff $p_i\Vdash \check{\delta}\in\name{a}$. The reason for defining $M$ is that we can guess its initial segments by our diamond sequence.
Let $S\subseteq S^{\kappa^+}_{\omega}$ be a stationary set for which $\eta\in S \Rightarrow M_\eta = M\cap\eta$, where $\langle M_\alpha: \alpha\in S^{\kappa^+}_{\omega} \rangle$ is the sequence of models derived from the diamond sequence that exemplifies $\Diamond_{S^{\kappa^+}_{\omega}}$.

Fix an ordinal $\eta\in S$, and recall that $T_\eta = \{\delta^\eta_n: n\in\omega\}$. By the choice of these ordinals, $\eta = \bigcup_{n\in\omega} \delta^\eta_n$. Likewise, $M\models P(\gamma^\eta_n,\gamma^\eta_{n+1}), R(\gamma^\eta_n,\delta^\eta_n)$ for every $n\in\omega$.
By the $M$-interpretation of $P$ and $R$ we can choose an increasing sequence $\langle q_n: n\in\omega\rangle$ of conditions in $\mathbb{P}$ such that $q_n \Vdash \delta^\eta_n\in\name{a}$, simply by letting $q_n = p_{\gamma^\eta_n}$ for every $n\in\omega$.

Every condition $q_n$ is an element of $N$, and by elementarity we see that $\langle q_n: n\in\omega\rangle$ is an increasing sequence also in $\mathcal{H}(\chi)$. Since $\mathbb{P}$ is $\aleph_1$-complete, there exists $q\in\mathbb{P}$ so that $\forall n\in\omega, q_n\leq q$. Notice that $p\leq q$ as all the $q_n$-s are above $p$.
It follows that $q\Vdash T_\eta\subseteq\name{a}$, and since $p$ was arbitrary this holds in the generic extension by $\mathbb{P}$, so we are done.

We use now the generalized Mathias forcing in order to increase $\mathfrak{s}_\kappa$ above $\kappa^+$.
This forcing notion is $\kappa$-directed-closed, so a fortiori $\aleph_1$-complete, and the value of $\mathfrak{s}_\kappa$ can be arbitrarily large.
In the generic extension we have $\clubsuit_{\kappa^+}$, since our tiltan is indestructible by the above forcing notion.
Likewise, $\mathfrak{s}_\kappa>\kappa^+$ has been forced, so the proof is accomplished.

\hfill \qedref{thmsupercompact}

The set $S^{\kappa^+}_{\omega}$ can be replaced, in the above proof, by $S^{\kappa^+}_{\theta}$ for every $\theta=\cf(\theta)<\kappa$.
Similarly, one can preserve tiltan at $S^\mu_\omega$ or $S^\mu_\theta$ for some $\theta=\cf(\theta)<\kappa$ when $\mu=\cf(\mu)>\kappa^+$, and force $\mathfrak{s}_\kappa\geq\mu$.

In the case of $\kappa=\aleph_0$ there is no infinite $\theta<\kappa$ for which we can build the tiltan over $S^{\aleph_1}_{\theta}$.
Hence we must step up to $\aleph_2$, work with $S^{\aleph_2}_{\aleph_0}$ and then collapse $\aleph_1$. This method is limited with respect to the splitting number.
It is also limited with respect to a very similar open problem, namely whether $\clubsuit$ is consistent with $\binom{\omega_1}{\omega} \rightarrow \binom{\omega_1}{\omega}^{1,1}_2$, see Question 3.21 in \cite{MR3509813}. Although it is opaque with respect to $\aleph_1$, we can answer positively at the supercompact level.

\begin{corollary}
\label{corpolarized} If $\kappa$ is supercompact then $\clubsuit_{\kappa^+}$ is consistent with the relation $\binom{\kappa^+}{\kappa} \rightarrow \binom{\kappa^+}{\kappa}^{1,1}_2$.
\end{corollary}

\par\noindent\emph{Proof}. \newline
Since $\mathfrak{s}_\kappa>\kappa^+$ implies $\binom{\kappa^+}{\kappa} \rightarrow \binom{\kappa^+}{\kappa}^{1,1}_2$, we can use the above theorem.

\hfill \qedref{corpolarized}

\newpage

\bibliographystyle{alpha}
\bibliography{arlist}

\end{document}